\documentclass[arxiv, preprint]{imsart}

 \RequirePackage{amsthm,amsmath,amsfonts,amssymb}
\RequirePackage[authoryear]{natbib} 
\RequirePackage[colorlinks,citecolor=blue,urlcolor=blue]{hyperref}
\RequirePackage{graphicx}
\usepackage{wrapfig}
\usepackage{multicol}
\usepackage{subcaption}
\usepackage{diagbox}

\startlocaldefs
\theoremstyle{plain}

\newtheorem{theorem}{Theorem}[section]

\theoremstyle{remark}

\newcommand{\ep}{\textsf{E}}

\endlocaldefs

\begin{document}

\begin{frontmatter}
\title{A note on the
permutation distribution of generalized correlation
coefficients}
\runtitle{Permutation distribution of generalized correlation
coefficients}

\begin{aug}
\author{\fnms{Yejiong} \snm{Zhu}\ead[label=e1]{yjzhu@ucdavis.edu}}
\and
\author{\fnms{Hao} \snm{Chen}\ead[label=e2]{hxchen@ucdavis.edu}}
\runauthor{Y. Zhu and H. Chen}
\address{University of California, Davis}
\end{aug}

\begin{abstract}
We provide sufficient conditions for the asymptotic normality of the generalized correlation coefficient $\sum a_{ij}b_{ij}$ under the permutation null distribution when $a_{ij}$'s are symmetric and $b_{ij}$'s are symmetric.
\end{abstract}



\end{frontmatter}


\section{Introduction}

The form of the generalized correlation coefficients $\sum a_{ij}b_{ij}$ is common under a variety of settings.  For example, given a sample $(x_i, y_i),i = 1,\cdots,N$, the popular Pearson and Spearman measures, ignoring standardization, can be obtained with $a_{ij} = x_i-x_j$, $b_{ij} = y_i-y_j$, and $a_{ij} = \text{rank}(x_i)-\text{rank}(x_j)$, $b_{ij} = \text{rank}(y_i)-\text{rank}(y_j)$, respectively.  It is also a unified form for non-parametric two-sample test when the permutation null distribution is of interest.  Explicitly, let $x_i$'s be the observations, $y_i$'s be 0's and 1's indicating the group ID, $m$ be the number of 0's in $y_i$'s, and $n$ be the number of 1's in $y_i$'s.  We give a few examples of non-parametric test statistics below.
\begin{enumerate}
\item Wilcoxon rank-sum/Mann-Whitney U test: $a_{ij} = \text{sign}(x_i-x_j)$, $b_{ij} = \text{sign}(y_i - y_j)$.
\item The edge-count test proposed in \cite{friedman1979}: Let $G$ be the minimum spanning tree constructed on $x_i$'s, then $a_{ij} = 1$ if $x_i$ and $x_j$ is connected in $G$, $b_{ij} = \text{abs}(y_i - y_j)$.
\item The MMD test proposed in \cite{gretton2012kernel}: Let $k(\cdot, \cdot)$ be a kernel function, then $a_{ij} = k(x_i, x_j)$, $b_{ij} =\frac{1}{m(m-1)}(1-y_i)(1-y_j) + \frac{1}{n(n-1)}y_iy_j - \frac{1}{mn}((1-y_i)y_j + y_i(1-y_j))$.
\item The weighted edge-count test proposed in \cite{chen2018weighted}: Let $G$ be a similarity graph constructed on $x_i$'s, then $a_{ij} = 1$ if $x_i$ and $x_j$ is connected in $G$,  $b_{ij} = (1-p)(1-y_i)(1-y_j) + py_iy_j$ with $p=m/N$ or $p=(m-1)/(N-2)$.
\end{enumerate}

Under the permutation null distribution, the test statistic becomes
\begin{align*}
    \Gamma = \sum_{i=1}^N\sum_{j=1}^Na_{ij}b_{\pi(i)\pi(j)},
\end{align*}
where $\pi(\dot)$ denotes a permutation of indices and each is with probability $\frac{1}{N!}$.  It is of generic interest in understanding the asymptotic distribution of $\Gamma$ as $N\rightarrow \infty$. 


\cite{daniels1944relation} provided sufficient conditions for studying the limiting distribution of $\Gamma$.  Their theorem is provided below.

\begin{theorem}{\citep{daniels1944relation}} \label{th1} 
Under conditions
\begin{align*}
    &a_{ij} =-a_{ji},\quad b_{ij} = -b_{ji},\quad a_{ii}=b_{ii} = 0,\\
    &\sum_{i,j,k}a_{ij}a_{ik} \asymp N^3(\max|a_{ij}|)^2,\quad \sum_{i,j,k}b_{ij}b_{ik} \asymp N^3(\max|b_{ij}|)^2,
\end{align*}
we have
\begin{align*}
    \frac{\Gamma}{2\sqrt{\frac{(\sum a_{ij}a_{ik})(\sum b_{ij}b_{ik})}{N^3}}}\rightarrow N(0,1),\text{ as } N \rightarrow \infty.
\end{align*}
\end{theorem}

The conditions in \cite{daniels1944relation} requires that $a_{ij} = -a_{ij}$, $b_{ij} = -b_{ij}$.  This is exactly the case for the Wilcoxon rank-sum test.  But for the other examples provided above, the asymmetric condition no longer holds.  One might do extra steps to make them asymmetric, while it would be more natural to remove the asymmetric condition.  In light of this, \cite{friedman1983graph} proposed a much relaxed condition of $\sum a_{ij} = \sum {b_{ij}} = 0$ in replacing the asymmetric condition.  They also relaxed the other conditions in \cite{daniels1944relation} to $\frac{\left(\sum_{i j k l} a_{i j} a_{i k} a_{i l}\right)^{2}}{\left(\sum_{i j k} a_{i j} a_{i k}\right)^{3}}\rightarrow0$ and $\frac{\left(\sum_{i j k l} b_{i j} b_{i k} b_{i l}\right)^{2}}{\left(\sum_{i j k} b_{i j} b_{i k}\right)^{3}}=0$.   However, they didn't provide a proof for their claims. 

Later, \cite{pham1989asymptotic} showed that the conditions in \cite{friedman1983graph} are not sufficient and they provide sufficient conditions for much more relaxed situations.  In particular, \cite{pham1989asymptotic} provided three set of sufficient conditions.  First, they define 
\begin{align*}
    a_{i j}^{\prime}&=\left(a_{i j}-\tfrac{1}{N(N-1)} \Sigma a_{k l}\right) I_{i\neq j},\\
    a_{i+}^{\prime}&=\sum_{j} a_{i j}^{\prime}, \quad a_{+j}^{\prime}=\sum_{i} a_{i j}^{\prime}, \\
   a_{i j}^{*}&=\left\{ \begin{array}{ll} \left(a_{i j}^{\prime}-\frac{1}{N-2}\left(a_{i+}^{\prime}+a_{+ j}^{\prime}\right)\right) I_{i\neq j}, & \text{ when } a_{ij} = a_{ji}, b_{ij} = b_{ji}, \\ \left( a_{i j}-\frac{1}{N}\left(a_{i+}^{\prime}+a_{+j}^{\prime}\right) \right) I_{i\neq j}, & \text{ when } a_{ij} = -a_{ji}, b_{ij} = -b_{ji},  \end{array} \right. 
\end{align*}
and $b_{ij}^\prime,b_{i+}^\prime,b_{+j}^\prime$ and $b_{ij}^*$ are defined similarly.

\allowdisplaybreaks
\begin{theorem}{\citep{pham1989asymptotic}}\label{Pham1}
When $a_{ij} = a_{ji}, b_{ij} = b_{ji}$ or $a_{ij} = -a_{ji}, b_{ij} = -b_{ji}$, under conditions 
\begin{align*}
    &a_{ii}= b_{ii} = 0, \forall i, \quad
    \frac{N\left(\Sigma a_{i j}^{* 2}\right)\left(\Sigma b_{i j}^{* 2}\right)}{\left(\sum a_{i+}^{\prime 2}\right)\left(\Sigma b_{i+}^{\prime 2}\right)} = o(1)\\
    &\frac{\sum_{i=1}^N (a_{i+}^\prime)^r}{\left(\sum_{i=1}^N (a_{i+}^\prime)^2\right)^\frac{r}{2}} = O(N^{1-\frac{r}{2}}), \text{ for all integer }r> 2\\
    &\frac{\sum_{i=1}^{N}\left|b_{i+}^\prime\right|^{r}}{\left[\sum_{i=1}^{N}\left(b_{i+}^\prime\right)^{2}\right]^{r/ 2}}=o(1),\text{ for some } r>2, \text{ or } \frac{\max(b_{i+}^\prime)^2}{\sum_{i=1}^N(b_{i+}^\prime)^2} = o(1),
\end{align*}
we have $\sum a_{ij}b_{\pi(i)\pi(j)}$ converges to a normal distribution as $N\rightarrow\infty$.
\end{theorem}

\begin{theorem}{\cite{pham1989asymptotic}} \label{Pham2}
When $a_{ij} = a_{ji}, b_{ij} = b_{ji}$ or $a_{ij} = -a_{ji}, b_{ij} = -b_{ji}$, under conditions 
\begin{align*}
    &\max \left(\sum_{j}\left|a_{i j}\right|\right)=O\left[\left(\sum a_{i j}^{2}\right) /\left(N \max \left|a_{i j}\right|\right)\right]\\
    &\left(\Sigma a_{i j}^{2}\right) /\left(N^{2} \max a_{i j}^{2}\right) = o(1),\\
    &\sum_{j} b_{i j}=0 \quad \text { for all } i,\quad \sum_{i} b_{i j}=0 \quad \text { for all } j,\\
    &\Sigma\left|b_{i j}\right|^{r} / N^{2}=O\left[\left(\Sigma b_{i j}^{2} / N^{2}\right)^{r / 2}\right], \quad \text{for all integers }r\geq 3,
\end{align*}
then $$
\frac{\sum a_{i j} b_{\pi(i) \pi(j)}}{\sqrt{2\left(\sum a_{i j}^{2}\right)\left(\sum b_{i j}^{2}\right) / N^{2}}} \rightarrow_d N(0,1), \text{ as } N\rightarrow \infty.$$
\end{theorem}

\begin{theorem}{\cite{pham1989asymptotic}} \label{Pham3}
When $a_{ij} = a_{ji}, b_{ij} = b_{ji}$ or $a_{ij} = -a_{ji}, b_{ij} = -b_{ji}$, under conditions
\begin{align*}
    &\sum a_{ij} = \sum b_{ij} = 0,\\
    &\max \left(\sum_{j}\left|a_{i j}\right|\right)=O\left(\max \left|a_{i j}\right|\right),\\
    &\liminf \left(\Sigma a_{i j}^{2}\right) /\left(N \max a_{i j}^{2}\right)>0,\\
    &\limsup \left(\sum b_{i+}^{2} / N\right) /\left(\Sigma b_{i j}^{2}\right)<1,
\end{align*}
then 
\begin{align*}
    \frac{\sum a_{ij}b_{\pi(i)\pi(j)}}{\sqrt{4\left(\sum^{\prime} a_{i j} a_{i k}\right)\left(\Sigma^{\prime} b_{i j} b_{i k}\right) / N^{3}+2\left(\Sigma a_{i j}^{2}\right)\left(\Sigma b_{i j}^{2}\right) / N^{2}}} \rightarrow_d N(0,1), \text{ as } N\rightarrow \infty.
\end{align*}
\end{theorem}

For the three set of conditions in \cite{pham1989asymptotic}, we can see that the conditions in Theorem \ref{Pham1} requires substantial computations to obtain $a_{ij}^*$'s  and $b_{ij}^*$'s before checking the conditions, which might be error-prone for the general community who might want to use the theorem.  On the other hand, the conditions in Theorem \ref{Pham2} and \ref{Pham3} are more user friendly, while these conditions are rather strong.  For example, if $a_{ij} = O(1), \forall i, j$, then $\sum a_{ij}^2 = O(N^2)$, $\max a_{ij}^2 = O(1)$, and the condition in Theorem \ref{Pham2} $\left(\Sigma a_{i j}^{2}\right) /\left(N^{2} \max a_{i j}^{2}\right) = o(1)$ would be violated.  Similarly, if $a_{ij} = O(1), \forall i, j$, then $\max (\sum_j |a_{ij}|) = O(N)$, $\max |a_{ij}| = O(1)$, and the condition in Theorem \ref{Pham3} $\max \left(\sum_{j}\left|a_{i j}\right|\right)=O\left(\max \left|a_{i j}\right|\right)$ would be violated.

In this paper, we provide a set of sufficient conditions for the symmetric setting ($a_{ij} = a_{ji}, b_{ij} = b_{ji}$).  The proof of the theorem can be extended from that in \cite{daniels1944relation}.  We write this note as the symmetric setting is common in many of those non-parametric test for high-dimensional data and non-Euclidean data \citep{friedman1979, gretton2012kernel, chen2017new, chen2018weighted}, and the set of conditions in our theorem is easy to use.

\section{Main theorem}

\begin{theorem} \label{main}
Under conditions
\begin{align*}
    &a_{ij} =a_{ji}, b_{ij} =b_{ji}, \sum a_{ij} = \sum{b_{ij}} = 0,\\
    &\sum_{i,j,k}a_{ij}a_{ik} \asymp N^3 (\max|a_{ij}|)^2,\quad \sum_{i,j,k}b_{ij}b_{ik} \asymp N^3(\max|b_{ij}|)^2,
\end{align*}
we have
\begin{align*}
    \frac{\Gamma}{2\sqrt{\frac{(\sum a_{ij}a_{ik})(\sum b_{ij}b_{ik})}{N^3}}}\rightarrow N(0,1),\text{ as } N \rightarrow \infty.
\end{align*}
\end{theorem}

We first check the scenario of $a_{ij} = O(1), \forall i, j$ that the conditions in Theorems \ref{Pham2} and \ref{Pham3} are violated.  In our case, $\sum_{i,j,k} a_{ij} a_{ik} = O(N^3), \max|a_{ij}| = O(1)$, then the condition $\sum_{i,j,k}a_{ij}a_{ik} \asymp N^3 (\max|a_{ij}|)^2$ holds.

\begin{proof}

Referring to the proof in \cite{daniels1944relation}, we firstly rewritten $\Gamma$ as 
\begin{align*}
    \Gamma = \sum_{i,j,k,l}p_{ij}p_{lk}a_{jk}b_{il},
\end{align*}
where $[p_{ij}]_{N*N}$ is a permutation matrix whose each row and column has only one nonzero entry $1$.

Then, we investigate the moments of $\Gamma$ under all permutations by starting from the first moment. The computation of first moment requires to count the number of repeats of the term $a_{jk}b_{il}$ with $p_{ij}=p_{lk}=1$ among all permutations. Note that if $i=l$, then the corresponding $j$ must be equal to $k$ as each row and column only has one nonzero entry. When $i\neq l$, the term $a_{jk}b_{il}$ repeats $(N-2)!$ over all permutation, and when $i = l$, the term $a_{jj}b_{ii}$ repeats $(N-1)!$ over permutations. Thus, the moment is 
\begin{align*}
    \ep(\Gamma) &= \frac{1}{N!}\left((N-2)!\sum_{i\neq l}\sum_{j\neq k}a_{jk}b_{il}+(N-1)!\sum a_{jj} \sum b_{ii}\right)\\
    &= \frac{1}{N-1}\sum b_{ii} \sum a_{ii}.
\end{align*}

For the second moment $\ep(\Gamma^2)$ , note that $(\sum_{i,j,k,l}p_{ij}p_{lk}a_{jk}b_{il})^2$ is the sum of $p_{ij}p_{lk}p_{rs}p_{tu}a_{jk}a_{su}b_{il}b_{rt}$. It can be computed by following the similar procedure by counting the number of repeats of $a_{jk}a_{su}b_{il}b_{rt}$ with $p_{ij}p_{lk}p_{rs}p_{tu}=1$ over all permutations. We call there $m$ independent subscripts if $m$ subscripts are not equal with each other. For example, if $i\neq l \neq r\neq t$, we call there are 4 independent subscribes. In addition, we define $\sum'$ as the sum with unequal subscripts, e.g. $\sum'a_{jk}a_{su} = \sum_{j}\sum_{k\neq j}\sum_{s\neq k,
s\neq j}\sum_{u\neq j, u \neq k, u\neq s} a_{jk}a_{su}$.

When there are 4 independent subscripts, the term $a_{jk}a_{su}b_{il}b_{rt}$ repeats $(N-4)!$ over all permutation. When there are 3 independent subscripts, the term $a_{jk}a_{ju}b_{il}b_{it}$ repeats $4(N-3)!$ and the term $a_{jj}a_{su}b_{ii}b_{rt}$ repeats $2(N-3)!$, where the number $4$ comes from the fact that the form $a_{jk}a_{ju}$ can be obtained by arranging subscripts in 4 possible ways since $a_{ij} = a_{ji}$, and the number 3 is due to the same reason. When there are 2 independent subscripts, the term $a_{jk}^2b_{il}^2$, $a_{jj}a_{ss}b_{ii}b_{rr}$ and $a_{jj}a_{ju}b_{ii}b_{it}$ repeats $2(N-2)!$, $(N-2)!$ and $4(N-2)!$ times, respectively. When there is only 1 independent subscript, the term $a_{jj}b_{ii}$ repeats $(N-1)!$ times. Thus, the second moment is 
\begin{align*}
    \ep(\Gamma^2) =& \frac{(N-4)!}{N!}(\sum{'}a_{jk}a_{su})(\sum{'}b_{il}b_{rt})+\frac{4(N-3)!}{N!}(\sum{'}a_{jk}a_{ju})(\sum{'}b_{il}b_{it})\\
    &+\frac{2(N-3)!}{N!}(\sum{'}a_{jj}a_{su})(\sum{'}b_{ii}b_{rt})+\frac{2(N-2)!}{N!}(\sum{'}a_{jk}^2)(\sum{'}b_{il}^2)\\
    &+\frac{(N-2)!}{N!}(\sum{'}a_{jj}a_{ss})(\sum{'}b_{ii}b_{rr})+\frac{4(N-2)!}{N!}(\sum{'}a_{jj}a_{ju})(\sum{'}b_{ii}b_{it})\\
    &+\frac{(N-1)}{N!}(\sum{'}a_{jj}^2)(\sum{'}b_{ii}^2)
\end{align*}

The $p$-th moment of $\Gamma$ can be computed by counting the number of repeats of the term $a_{j_1k_1}\cdots a_{j_pk_p}*b_{i_1l_1}\cdots b_{i_pl_p}$ under different cases. It is not hard to see that it is the sum of 
\begin{align}
   \frac{(N-f)!}{N!}A_f \sum{'}a_{j_1k_1}\cdots a_{j_pk_p}\sum{'}b_{i_1l_1}\cdots b_{i_pl_p} 
   \label{term1}
\end{align}
where $f$ is the number of independent subscripts in the $\sum{'}$ and $A_f$ is the number of possible arrangements to obtain a specific form of $a_{j_1k_1}\cdots a_{j_pk_p}$, which does not depend on $N$.


Let $a_{max} = \max|a_{ij}|$ and $b_{max} = \max|b_{ij}|$.
Note that the sum $\sum{'}$ can be reorganized as the linear combination of the corresponding sum $\sum$ with the same subscripts and other $\sum$'s with additional tied subscripts, e.g. $\sum{'}a_{jk} = \sum a_{jk}-\sum a_{jj}$. In the sum $\sum$, if one $a_{jk}$ contains unique subscripts $j,k$ that do not repeat in other $a_{..}$'s, the sum $\sum$ would be zero because of the assumption $\sum a_{ij} = 0$.

For even moments, i.e. $p=2m$, \cite{daniels1944relation} showed that the sum $\sum a_{j_1k_1}\cdots a_{j_{2m}k_{2m}}$ would not vanish with at most $3m$ independent subscripts. In the symmetric case $a_{ij} = a_{ji}$, the sum with $3m$ independent subscripts can always be arranged as 
\begin{align*}
    \sum a_{j_1k_1}a_{j_1k_2}\cdots a_{j_mk_{2m-1}}a_{j_mk_{2m}} = (\sum a_{j_1k_1}a_{j_1k_2})^m,
\end{align*}
which is of the order of $N^{3m}a_{max}^{2m}$ under the condition $\sum_{i,j,k}a_{ij}a_{ik} \asymp N^3a_{max}^2$. Any other sum $\sum a_{j_1k_1}\cdots a_{j_{2m}k_{2m}}$ with at most $3m-1$ independent subscripts is of the order at most $N^{3m-1}a_{max}^{2m}$ as there are fewer than $3m$ summations from 1 to $N$. Similar results also hold for the sum of $b_{il}$'s. Hence, the term (\ref{term1}) with $f = 3m$ is of the order $$N^{-3m}N^{3m}N^{3m}a_{max}^{2m}b_{max}^{2m} = N^{3m}a_{max}^{2m}b_{max}^{2m}.$$
When $f\leq 3m-1$, the order of the term (\ref{term1}) is at most $N^{-f}N^{f}N^{f}a_{max}^{2m}b_{max}^{2m}$ that is dominated by $N^{3m}a_{max}^{2m}b_{max}^{2m}$. When $f\geq 3m+1$, the order of the term (\ref{term1}) is at most $N^{-f}N^{3m}a_{max}^2N^{3m}b_{max}^2$ that is also dominated by $N^{3m}a_{max}^{2m}b_{max}^{2m}$. Thus, the $2m$-th moment is of the order of $N^{3m}a_{max}^2b_{max}^2$.

For odd moments, i.e. $p= 2m+1$, \cite{daniels1944relation} showed that the sum $\sum a_{j_1k_1}\cdots a_{j_{2m+1}k_{2m+1}}$ would not vanish with at most $3m+1$ independent subscripts and the order of the sum with $3m+1$ independent subscripts is at most $N^{3m+1}a_{max}^{2m+1}$. The sum consisting of $b_{il}$ has the similar result. Thus, the term (\ref{term1}) with $f = 3m+1$ is of the order of at most $N^{3m+1}a_{max}^{2m+1}b_{max}^{2m+1}$. When $f\leq 3m$, the order the term (\ref{term1}) is at most $N^{-f}N^fa_{max}^{2m+1}N^{f}b_{max}^{2m+1}$ that is dominated by $N^{3m+1}a_{max}^{2m+1}b_{max}^{2m+1}$. When $f\geq 3m+2$, the term (\ref{term1}) is of the order of at most $N^{-f}N^{3m+1}a_{max}^{2m+1}N^{3m+1}b_{max}^{2m+1}$ that is also dominated by $N^{3m+1}a_{max}^{2m+1}b_{max}^{2m+1}$. Thus, $2m+1$-th moment is of the order of at most $N^{3m+1}a_{max}^{2m+1}b_{max}^{2m+1}$.

Consider the scaled version $\frac{\Gamma}{N^{1.5}a_{max}b_{max}}$, the $(2m+1)$-th term in the Taylor's expansion of the moment generating function $\ep e^{t\Gamma N^{-\frac{3}{2}}a_{max}^{-1}b_{max}^{-1}}$ is
\begin{align*}
    \frac{t^{2m+1}}{(2m+1)!}\ep(\Gamma^{2m+1})N^{-3m-\frac{3}{2}}a_{max}^{-2m-1}b_{max}^{-2m-1} \precsim \frac{t^{2m+1}}{(2m+1)!}N^{-\frac{1}{2}},
\end{align*}
which goes to zero when $N$ goes to infinity. The $2m$-th term in the Taylor's expansion of the moment generating function $\ep e^{t\Gamma N^{-\frac{3}{2}}a_{max}^{-1}b_{max}^{-1}}$ is
\begin{align*}
    \frac{t^{2m}}{(2m)!}\ep(\Gamma^{2m})N^{-3m}a_{max}^{-2m}b_{max}^{-2m} &\asymp \frac{t^{2m}}{(2m)!}A_{2m}\left(\frac{\sum a_{j_{1} k_{1}} a_{j_{1} k_{2}}}{N^{3}a_{max}^2}\right)^{m}\left(\frac{\sum b_{i_{1} l_{1}} b_{i_{1} l_{2}}}{N^3b_{max}^2}\right)^{m}\\
    &= \frac{(2t^2)^m}{m!}\left(\frac{\sum a_{j_{1} k_{1}} a_{j_{1} k_{2}}}{N^{3}a_{max}^2}\right)^{m}\left(\frac{\sum b_{i_{1} l_{1}} b_{i_{1} l_{2}}}{N^3b_{max}^2}\right)^{m},
\end{align*}
where the last equality is due to $A_{2m} = \frac{(2m)!2^{m}}{m!}$ from \cite{daniels1944relation}. Thus, we have 
\begin{align*}
    \frac{\Gamma}{N^{1.5}a_{max}b_{max}} \rightarrow_d N(0, 4h_ah_b),\text{ as } N\rightarrow \infty
\end{align*}
with $h_a = \lim_{N\rightarrow \infty} \frac{\sum a_{ij} a_{ik}}{N^{3}a_{max}^2}$ and $h_b = \lim_{N\rightarrow \infty} \frac{\sum b_{ij} b_{ik}}{N^{3}b_{max}^2}$, and it is equivalent as
\begin{align*}
    \frac{\Gamma}{2\sqrt{\frac{(\sum a_{ij}a_{ik})(\sum b_{ij}b_{ik})}{N^3}}}\rightarrow_d N(0,1),\text{ as } N \rightarrow \infty.
\end{align*}

\end{proof}

\bibliographystyle{imsart-nameyear} 
\bibliography{reference}      


\end{document}